\magnification1200
\input amssym.def
\input amssym.tex

\hsize=13.5truecm
\baselineskip=16truept plus.3pt minus .3pt

\font\secbf=cmb10 scaled 1200
\font\eightrm=cmr8
\font\sixrm=cmr6

\font\eighti=cmmi8

\font\sixi=cmmi6
\skewchar\eighti='177 \skewchar\sixi='177

\font\eightsy=cmsy8
\font\sixsy=cmsy6
\skewchar\eightsy='60 \skewchar\sixsy='60

\font\eightit=cmti8

\font\eightbf=cmbx8
\font\sixbf=cmbx6

\let\sc=\tensc

\font\eightsc=cmcsc10 scaled 800
\font\secbf=cmb10 scaled 1200
\font\subsecfont=cmb10 scaled \magstephalf
\font\amb=cmmib10

\font\ambi=cmmib10 scaled 700

\newfam\mbfam 

\textfont\mbfam\amb \scriptfont\mbfam\ambi


\def\aa{\def\rm{\fam0\eightrm}%
  \textfont0=\eightrm \scriptfont0=\sixrm \scriptscriptfont0=\fiverm
  \textfont1=\eighti \scriptfont1=\sixi \scriptscriptfont1=\fivei
  \textfont2=\eightsy \scriptfont2=\sixsy \scriptscriptfont2=\fivesy
  \textfont3=\tenex \scriptfont3=\tenex \scriptscriptfont3=\tenex
  \def\sc{\eightsc}
  \def\it{\fam\itfam\eightit}%
  \textfont\itfam=\eightit
  \def\bf{\fam\bffam\eightbf}%
  \textfont\bffam=\eightbf \scriptfont\bffam=\sixbf
   \scriptscriptfont\bffam=\fivebf
  \normalbaselineskip=9.7pt
  \setbox\strutbox=\hbox{\vrule height7pt depth2.6pt width0pt}%
  \normalbaselines\rm}

\def\Proof{\vskip12pt\noindent{\bf Proof.} }

\def\Remark#1{\vskip12pt\noindent{\bf Remark #1}}

\def\m@th{\mathsurround=0pt}

\def\cc#1{\hbox to .89\hsize{$\displaystyle\hfil{#1}\hfil$}\cr}
\def\lc#1{\hbox to .89\hsize{$\displaystyle{#1}\hfill$}\cr}
\def\rc#1{\hbox to .89\hsize{$\displaystyle\hfill{#1}$}\cr}

\def\eqal#1{\null\,\vcenter{\openup\jot\m@th
  \ialign{\strut\hfil$\displaystyle{##}$&&$\displaystyle{{}##}$\hfil
      \crcr#1\crcr}}\,}

\def\section#1{\vskip 22pt plus6pt minus2pt\penalty-400
        {{\secbf
        \noindent#1\rightskip=0pt plus 1fill\par}}
        \par\vskip 12pt plus5pt minus 2pt
        \penalty 1000}

\def\subsection#1{\vskip 20pt plus6pt minus2pt\penalty-400
        {{\subsecfont
        \noindent#1\rightskip=0pt plus 1fill\par}}
        \par\vskip 8pt plus5pt minus 2pt
        \penalty 1000}

\def\subsubsection#1{\vskip 18pt plus6pt minus2pt\penalty-400
        {{\subsecfont
        \noindent#1}}
        \par\vskip 7pt plus5pt minus 2pt
        \penalty 1000}

\def\centerlast#1{\begingroup \leftskip=0pt plus 1fil\rightskip=0pt plus -1fil
\parfillskip=0pt plus 2fil
\parindent 0pt
\par #1\par\endgroup}
\def\\{\hfill\break}

\def\kwadrat{\null\ \hfill\null\ \hfill$\square$}
\def\mida#1{{{\null\kern-4.2pt\left\bracevert\vbox to 6pt{}\!\hbox{$#1$}\!\right\bracevert\!\!}}}
\def\midy#1{{{\null\kern-4.2pt\left\bracevert\!\!\hbox{$\scriptstyle{#1}$}\!\!\right\bracevert\!\!}}}

\def\diagin{\hbox{---}\hskip-11.5pt\intop}
\def\diagint{\hbox{--}\hskip-8.2pt\intop}
\def\diagintop{\mathop{\mathchoice
{{\diagin}}%
{{\diagint}}%
{{\diagint}}%
{{\diagint}}%
}\limits}

\def\divv{{{\rm div}\,}}
\def\rot{{\rm rot}\,}

\def\today{${\scriptscriptstyle\number\day-\number\month-\number\year}$}
\footline={{\hfil\rm\the\pageno\hfil${\scriptscriptstyle\rm\jobname}$\ \ \today}}

\def\ifnextchar#1#2#3{\bgroup
  \def\reserveda{\ifx\reservedc #1 \aftergroup\firstoftwo
    \else \aftergroup\secondoftwo\fi\egroup{#2}{#3}}%
  \futurelet\reservedc\ifnch
  } 
\def\ifnch{\ifx \reservedc \sptoken \expandafter\xifnch
      \else \expandafter\reserveda
      \fi} 
\def\firstoftwo#1#2{#1}
\def\secondoftwo#1#2{#2} 
\def\tempswafalse{\let\iftempswa\iffalse}
\def\tempswatrue{\let\iftempswa\iftrue} 

\def\cite{\ifnextchar [{\tempswatrue\citea}{\tempswafalse\citeb}}
\def\citea[#1]#2{[#2, #1]}
\def\citeb#1{[#1]}\phantom{]}
\def\\{\hfil\break}

\def\N{{\Bbb N}}
\def\R{{\Bbb R}}
\def\T{{\Bbb T}}

\def\divv{{\rm div}\,}
\def\rot{{\rm rot}}

\centerlast{\secbf Some stability problem to the Navier-Stokes equations in 
the periodic case}

\vskip1cm
\centerline{\bf W. M. Zaj\c aczkowski}
\vskip1cm
\noindent Institute of Mathematics Polish Academy of Sciences,\\ 
\'Sniadeckich 8, 00-656 Warsaw, Poland,\\
e-mail:wz@impan.pl;\\
Institute of Mathematics and Cryptology, Cybernetics Faculty,\\ 
Military University of Technology,\\
Kaliskiego 2, 00-908 Warsaw, Poland
\vskip1cm

\noindent
{\bf Abstract.} 
The Navier-Stokes motions in a box with periodic boundary conditions are 
considered. First the existence of global regular two-dimensional solutions 
is proved. The solutions are such that continuous with respect to time norms 
are controlled by the same constant for all $t\in\R_+$. Assuming that the 
initial velocity and the external force are sufficiently close to the initial 
velocity and the external force of the two-dimensional solutions we prove 
existence of global three-dimensional regular solutions which remain close to 
the two-dimensional solutions for all time. In this way we mean stability of 
two-dimensional solutions. 


\vfil\eject 

\section{1. Introduction}

The aim of this paper is to prove stability of two-dimensional periodic solutions in the set of 
three-dimensional periodic solutions to the Navier-Stokes equation. We consider the three-dimensional fluid motions in the box $\Omega=[0,L]^3$, $L>0$, described by
$$\eqal{
&v_t+v\cdot\nabla v-\nu\Delta v+\nabla p=f\quad &{\rm in}\ \ \Omega\times\R_+,\cr 
&\divv v=0\quad &{\rm in}\ \ \Omega\times\R_+,\cr 
&v|_{t=0}=v(0)\quad &{\rm in}\ \ \Omega,\cr}
\leqno(1.1)
$$
where $v=(v_1(x,t),v_2(x,t),v_3(x,t))\in\R^3$ is the velocity of the\break 
fluid, $x=(x_1,x_2,x_3)$ with $x_i\in(0,L)$, $i=1,2,3$, is a given Cartesian\break 
system of coordinates, $p=p(x,t)\in\R$ is the pressure and $f=\break
(f_1(x,t),f_2(x,t),f_3(x,t))\in\R^3$ is the external force field. Finally, 
$\nu>0$ is the constant viscosity coefficient and the dot denotes the scalar 
product in $\R^3$.

\noindent
Two-dimensional solutions to (1.1) are such that $v=v_s=(v_{s1}(x_1,x_2,t)$, 
$v_{s2}(x_1,x_2,t),0)\in\R^2$, \ $p=p_s(x_1,x_2,t)\in\R$, \ $f=f_s=(f_{s1}(x_1,x_2,t),\break f_{s2}(x_1,x_2,t),0)\in\R^2$ and satisfy the problem
$$\eqal{
&v_{st}+v_s\cdot\nabla v_s-\nu\Delta v_s+\nabla_{p_s}=f_s\quad &{\rm in}\ \ \Omega\times\R_+,\cr 
&\divv v_s=0\quad &{\rm in}\ \ \Omega\times\R_+,\cr 
&v_s|_{t=0}=v_s(0)\quad &{\rm in}\ \ \Omega.\cr}
\leqno(1.2)
$$
To show stability, we introduce the quantities
$$
u=v-v_s,\quad q=p-p_s
\leqno(1.3)
$$
which are solutions to the problem
$$\eqal{
&u_t+u\cdot\nabla u-\nu\Delta u+\nabla q=-v_s\cdot\nabla u-u\cdot\nabla v_s+g\quad &{\rm in}\ \ \Omega\times\R_+,\cr 
&\divv u=0\quad &{\rm in}\ \ \Omega\times\R_+,\cr 
&u|_{t=0}=u(0)\quad &{\rm in}\ \ \Omega,\cr}
\leqno(1.4)
$$
with $g=f-f_s$. Therefore, to show stability of solutions to (1.2) we need to 
prove smallness of quantities (1.3) in some norms for all $t\in\R_+$. For this 
purpose we apply the energy method. For this we need the Poincar\'e 
inequality. Since it does not hold for solutions to problems (1.2) and (1.4) 
we introduce the quantities
$$\eqal{
&\bar v_s=v_s-\diagintop_\Omega v_sdx,\quad 
&\bar p_s=p_s-\diagintop_\Omega p_sdx,\quad 
&\bar f_s=f_s-\diagintop f_sdx,\cr 
&\bar u=u-\diagintop_\Omega udx,\quad &\bar q=q-\diagintop_\Omega qdx,\quad 
&\bar g=g-\diagintop_\Omega gdx,\cr}
\leqno(1.5)
$$
where
$$
\diagintop_\Omega\omega dx={1\over|\Omega|}\intop_\Omega\omega dx\quad 
{\rm and}\ \ |\Omega|=L^3.
$$
Applying the mean operator to $(1.2)_1$ and $(1.4)_1$ and using the periodic boundary conditions we have
$$
{d\over dt}\diagintop_\Omega v_sdx=\diagintop_\Omega f_sdx
\leqno(1.6)
$$
and
$$
{d\over dt}\diagintop_\Omega udx=\diagintop_\Omega gdx.
\leqno(1.7)
$$
In view of (1.6), (1.7) and that any space derivative of the mean vanishes we see that for quantities (1.5) problems (1.2) and (1.4) take the forms
$$\eqal{
&\bar v_{st}+v_s\cdot\nabla\bar v_s-\nu\Delta\bar v_s+\nabla\bar p_s=\bar f_s\quad &{\rm in}\ \ \Omega\times\R_+,\cr 
&\divv\bar v_s=0\quad &{\rm in}\ \ \Omega\times\R_+,\cr 
&\bar v_s|_{t=0}=\bar v_s(0)\quad &{\rm in}\ \ \Omega,\cr}
\leqno(1.8)
$$
and
$$\eqal{
&\bar u_t+u\cdot\nabla\bar u-\nu\Delta\bar u+\nabla\bar q=
-v_s\cdot\nabla\bar u-u\cdot\nabla\bar v_s+\bar g\quad &{\rm in}\ \ \Omega\times\R_+,\cr 
&\divv\bar u=0\quad &{\rm in}\ \ \Omega\times\R_+,\cr 
&\bar u|_{t=0}=\bar u(0)\quad &{\rm in}\ \ \Omega.\cr}
\leqno(1.9)
$$
Now, we formulate the main results of this paper (for any notation see 
Section 2). From Lemmas 3.1--3.3 we have

\proclaim Theorem 1. {\rm(two-dimensional solutions)}
Assume that $\bar f_s\in L_2(kT,\break(k+1)T; L_\sigma(\Omega))$, $k\in\N_0$, $\sigma>3$ and 
$\bar v_s(0)\in B_{\sigma,2}^1(\Omega)$. Then there exists a solution to problem (1.8) such that $\bar v_s\in W_{\sigma,2}^{2,1}(\Omega\times(kT,(k+1)T))$, 
$\nabla\bar p_s\in L_2(kT,(k+1)T;L_\sigma(\Omega))$, $k\in\N_0$ and the estimates hold
$$
\|\bar v_s\|_{W_{\sigma,2}^{2,1}(\Omega\times(0,T))}\le c(A^2+
\|\bar f_s\|_{L_2(0,T;L_\sigma(\Omega))}+\|v_s(0)\|_{B_{\sigma,2}^1(\Omega)})
\leqno(1.10)
$$
and
$$
\|\bar v_s\|_{W_{\sigma,2}^{2,1}(\Omega\times(kT,(k+1)T))}\le c\bigg({1\over\delta}A+A^2+
\|\bar f_s\|_{L_2((k-1)T,(k+1)T;L_\sigma(\Omega))}\bigg)
\leqno(1.11)
$$
where $k\in\N$, $\delta\in(1/2,1)$, 
$A=A_5=\big(1+{c_{s1}\over1-\exp(-\nu c_{s1}T)}\big)A_1^2+\|v_{sx}(0)\|_{L_2}^2$. 
Moreover, 
$A_1^2={1\over\nu c_{s1}}\sum_{k\in\N_0}\intop_{kT}^{(k+1)T}\|\bar f_s(t)\|_{L_2}^2dt$, $c_{s1}$ is the constant from the Poincar\'e inequality (2.3).
\goodbreak

\noindent
Lemma 4.2 implies

\proclaim Theorem 2. {\rm (stability)}
Let the assumptions of Theorem 1 hold. Let $\gamma\in(0,\gamma_*]$, where 
$\nu c_4-{c_5\over\nu^3}\gamma_*^2\ge{c_*\over2}$, $c_*<\nu c_4$
where $c_4(c_1)$ and $c_1$ is the constant from the Poincar\'e inequality 
(2.4). Assume that $\bar g\in C(\R_+;L_2)$, $\bar u(0)\in H^1$. 
Assume that
$$\eqal{
&\|\bar u(0)\|_{H^1}^2\le\gamma,\cr 
&G^2(t)\equiv{c_5\over\nu}\bigg[
\|\bar v_{sx}\|_{L_3}^2\bigg|\intop_0^t\diagintop_\Omega 
g(x,t')dxdt'+\diagintop_\Omega u(0)dx\bigg|^2
+\|\bar g\|_{L_2}^2\bigg]\le c_*{\gamma\over4},\cr} 
\leqno(1.13)
$$
Let $T>0$ be given and $k\in\N_0$. Assume that
$$\eqal{
&{c_5\over\nu}\intop_{kT}^{(k+1)T}\|\bar v_{sx}\|_{L_3}^2dt\le{c_*\over4}T,\qquad \intop_{kT}^{(k+1)T}G^2(t)dt\le\alpha\gamma,\cr 
&\alpha\exp\bigg({c_*\over4}T\bigg)\exp\bigg(-{c_*\over4}T\bigg)\le 1\cr}
$$
where constant $c_5$ appears in (4.16). Then
$$
\|\bar u(t)\|_{H^1}^2\le\gamma\quad {\rm for}\ \ t\in\R_+.
\leqno(1.14)
$$

\noindent
Finally by the regularity theory to the Navier-Stokes equations we have

\proclaim Theorem 3.  
Let the assumptions of Theorems 1 and 2 hold. Then there exists a solution to 
problem (1.1) such that $v=v_s+u\in W_2^{2,1}(\Omega\times(kT,\break(k+1)T))$,
$\nabla p=\nabla(p_s+q)\in L_2(kT,(k+1)T;L_2(\Omega))$, $k\in\N_0$, 
where $v_s$, $p_s$, $u$ are determined by Theorems 1 and 2, respectively.
\goodbreak

The first results connected with the stability of global regular solutions to 
the nonstationary Navier-Stokes equations were proved by Beirao da Veiga and 
Secchi \cite{1}, followed by Ponce, Racke, Sideris and Titi \cite{2}. 
Paper \cite{1} is concerned with the stability in $L_p$-norm of a strong 
three-dimensional solution of the Navier-Stokes system with zero external 
force in the whole space. In \cite{2}, assuming that the external force is 
zero and a three-dimensional initial function is close to 
a two-dimensional one in $H^1(\R^3)$, the authors showed the existence of 
a global strong solution in $\R^3$ which remains close to a two-dimensional 
strong solution for all times. In \cite{3} Mucha obtained a similar result 
under weaker assumptions about the smallness of the initial velocity 
perturbation.

In the class of weak Leray-Hopf solutions the first stability result was 
obtained by Gallagher \cite{4}. She proved the stability of 
two-dimensional solutions of the Navier-Stokes equations with 
periodic boundary conditions under three-dimensional perturbations both in 
$L_2$ and $H^{1\over2}$ norms.

The stability of nontrivial periodic regular solutions to the Navier-Stokes 
equations was studied by Iftimie \cite{5} and by Mucha \cite{6}. 
The paper \cite{6} is devoted to the case when the external force is 
a potential belonging to $L_{r,loc}(\T^3\times[0,\infty))$ and when the initial 
data belongs to the space $W_r^{2-2/r}(\T^3)\cap L_2(\T^3)$, where $r\ge2$ and 
$\T$ is a torus. Under the assumption that there exists a global solution with 
data of regularity mentioned above and that small perturbations 
of data have the same regularity as above, the author proves that 
perturbations of the velocity and the gradient of the pressure remain small in 
the spaces $W_r^{2,1}(\T^3\times(k,k+1))$ and $L_r(\T^3\times(k,k+1))$, 
$k\in\N$, respectively. Paper \cite{5} contains results concerning the 
stability of two-dimensional regular solutions to the Navier-Stokes system in 
a three-dimensional torus but here the initial data in the three-dimensional 
problem belongs to an anisotropic space of functions having different 
regularity in the first two directions than in the third direction, and the 
external force vanishes. Moreover, Mucha \cite{7} studies the stability of 
regular solutions to the nonstationary Navier-Stokes system in $\R^3$ assuming 
that they tend in $W_r^{2,1}$ spaces $(r\ge2)$ to constant flows.

The papers of Auscher, Dubois and Tchamitchian \cite{8} and of Gallagher, 
Iftimie and Planchon \cite{9} concern the stability of global regular 
solutions to the Navier-Stokes equations in the whole space $\R^3$ with zero 
external force. These authors assume that the norms of the considered 
solutions decay as $t\to\infty$.

It is worth mentioning the paper of Zhou \cite{10}, who proved the asymptotic 
stability of weak solutions $u$ with the property: 
$u\in L_2(0,\infty,BMO)$ to the Navier-Stokes equations in $\R^n$, $n\ge3$, 
with force vanishing as $t\to\infty$.

An interesting result was obtained by Karch and Pilarczyk \cite{11}, who 
concentrate on the stability of Landau solutions to the Navier-Stokes system 
in $\R^3$. Assuming that the external force is a singular distribution they 
prove the asymptotic stability of the solution under any $L_2$-perturbation.

Paper \cite{12} of Chemin and Gallagher is devoted to the stability 
of some unique global solution with large data in a very weak sense.

Finally, the stability of Leray-Hopf weak solutions has recently been examined 
by Bardos et al. \cite{13}, where equations with vanishing external force are 
considered. That paper concerns the following three cases: two-dimensional 
flows in infinite cylinders under three-dimensional perturbations which are 
periodic in the vertical direction; helical flows in circular cylinders under 
general three-dimensional perturbations; and axisymmetric flows under general 
three-dimensional perturbations. The theorem concerning the first case 
extends a result obtained by Gallagher \cite{4} for purely periodic 
boundary conditions.

\noindent
Most of the papers discussed above concern to the case with zero external 
force \cite{1, 2, 3, 5, 8, 9, 12, 13},
or with force which decays as $t\to\infty$ (\cite{10}). Exceptions are 
\cite{6, 7, 11}, where very special external forces, which are singular 
distributions in \cite{11} or potentials in \cite{6, 7}, are 
considered. However, the case of potential forces is easily reduced to the 
case of zero external forces.

The aim of our paper is to prove the stability result for a large class of 
external forces $f_s$ which do not produce solutions decaying as $t\to\infty$.

It is essential that our stability results are obtained together with the 
existence of a global strong three-dimensional solution close to 
a two-dimensional one.

The paper is divided into two main parts. In the first we prove existence of 
global strong two-dimensional solutions not vanishing as $t\to\infty$ because 
the external force does not vanish either. To prove existence of such 
solutions we use the step by step method. For this purpose we have to show 
that the data in the time interval $[kT,(k+1)T]$, $k\in\N$, do not increase 
with $k$. We do not need any restrictions on the time step $T$. 

\noindent
In the second part we prove existence of three-dimensional solutions that 
remain close to two-dimensional solutions. For this we need the initial 
velocity and the external force to be sufficiently close in apropriate norms 
to the initial velocity and the external force of the two-dimensional 
problems.

The proofs of this paper are based on the energy method, which strongly 
simplifies thanks to the periodic boundary conditions. The proofs of global 
existence which follow from the step by step technique are possible thanks to 
the natural decay property of the Navier-Stokes equations. This is mainly used 
in the first part of the paper (Section 3). To prove stability (Section 4) we 
use smallness of data $(v(0)-v_s(0)),(f-f_s)$ and a contradiction argument 
applied to the nonlinear ordinary differential inequality (4.20).

We restrict ourselves to proving estimates only, because existence follows 
easily by the Faedo-Galerkin method.

The paper is a serious generalization of \cite{14} because proofs are simpler, 
there is imposed less restrictions on data and there is no relation between 
$T$, $\nu$ and $f_s$ which in \cite{14} implies some smallness for 
two-dimensional solutions.

The paper is organized as follows. In Section 2 we introduce notation and give 
some auxiliary results. Section 3 is devoted to the existence of 
a two-dimensional solution. It also contains some useful estimates of the 
solution. In Section 4 we prove the existence of a global strong solution 
to problem (1.1) close to the two-dimensional solution for all time.

\section{2. Notation and auxiliary results}

By $L_p(\Omega)$, $p\in[1,\infty]$, we denote the Lebesgue space of integrable functions and by $H^s(\Omega)$, $s\in\N_0=\N\cup\{0\}$, the Sobolev space of functions with the finite norm
$$
\|u\|_{H^s}\equiv\|u\|_{H^s(\Omega)}=
\sum_{|\alpha|\le s}\bigg(\intop_\Omega|D_x^\alpha u|^2dx\bigg)^{1/2},
$$
where $D_x^\alpha=\partial_{x_1}^{\alpha_1}\partial_{x_2}^{\alpha_2}\partial_{x_3}^{\alpha_3}$, $|\alpha|=\alpha_1+\alpha_2+\alpha_3$, $\alpha_i\in\N_0$, $i=1,2,3$.

\proclaim Lemma 2.1. 
Assume that $\diagintop_\Omega f_s(t)dx$, $\diagintop_\Omega g(t)dx$ are locally integrable on $\R_+$ and $\diagintop_\Omega v_s(0)dx$, $\diagintop_\Omega u(0)dx$ are finite. Then for all $t\in\R_+$,
$$
\diagintop_\Omega v_s(t)dx=\intop_0^t\diagintop_\Omega f_s(t')dxdt'+\diagintop_\Omega v_s(0)dx,
\leqno(2.1)
$$
$$
\diagintop_\Omega u(t)dx=\intop_0^t\diagintop_\Omega g(t')dxdt'+\diagintop_\Omega u(0)dx.
\leqno(2.2)
$$

\Proof 
Applying the mean operator to (1.2) and (1.4), integrating by parts and using the periodic boundary conditions, we get (2.1) and (2.2) after integration with respect to time, respectively. This concludes the proof.

\proclaim Lemma 2.2. 
By the Poincar\'e inequality holds
$$
c_{s1}\|\bar v_s\|_{H^1}^2\le\|\nabla\bar v_s\|_{L_2}^2
\leqno(2.3)
$$
and
$$
c_1\|\bar u\|_{H^1}^2\le\|\nabla\bar u\|_{L_2}^2,
\leqno(2.4)
$$
where $c_{s1}$, $c_1$ are positive constants.

\noindent
Let as introduce the anisotropic Lebesgue and Sobolev spaces with the mixed norms, $L_{p_1,p_2}(\Omega\times(0,T))$ and $W_{p_1,p_2}^{2,1}(\Omega\times(0,T))$, $p_1,p_2\in(1,\infty)$, with the following norms
$$
\|u\|_{L_{p_2}(0,T;L_{p_1}(\Omega))}\equiv\|u\|_{L_{p_1,p_2}(\Omega\times(0,T))}=
\bigg(\intop_0^T\bigg(\intop_\Omega|u|^{p_1}dx\bigg)^{p_2/p_1}dt\bigg)^{1/p_2},
$$
$$\eqal{
&\|u\|_{W_{p_1,p_2}^{2,1}(\Omega\times(0,T))}=
\|D_x^2u\|_{L_{p_1,p_2}(\Omega\times(0,T))}+\|\partial_tu\|_{L_{p_1,p_2}(\Omega\times(0,T))}\cr 
&\quad+\|u\|_{L_{p_1,p_2}(\Omega\times(0,T))}.\cr}
$$
We introduce the Besov space $B_{p,q}^s(\Omega)$ (see [\cite{15}, Ch. 7, 
Sect. 7.32]) by
$$
B_{p,q}^s(\Omega)=(L_p(\Omega),W_p^m(\Omega))_{s/m,q,J}.
$$
In \cite[Ch. 4, Sect. 18] {16} the Besov spaces are introduced more explicitly.\\
Let us consider the Stokes system
$$\eqal{
&\omega_t-\nu\Delta\omega+\nabla q=f\quad &{\sl in}\ \ \Omega\times(0,T),\cr 
&\divv\omega=0\quad &{\sl in}\ \ \Omega\times(0,T),\cr 
&\omega|_{t=0}=\omega(0)\quad &{\sl in}\ \ \Omega.\cr}
\leqno(2.5)
$$
\kwadrat

\proclaim Lemma 2.3. 
Let $f\in L_{p_2}(0,T;L_{p_1}(\Omega))$, $p_1,p_2\in(1,\infty)$, \\
$\omega(0)\in B_{p_1,p_2}^{2-2/p_2}(\Omega)$. Then there exists a solution to 
problem (2.5) such that $\omega\in W_{p_1,p_2}^{2,1}(\Omega\times(0,T))$, 
$\nabla q\in L_{p_2}(0,T;L_{p_1}(\Omega))$ and
$$\eqal{
&\|\omega\|_{W_{p_1,p_2}^{2,1}(\Omega\times(0,T))}+
\|\nabla q\|_{L_{p_2}(0,T;L_{p_1}(\Omega))}\le c
(\|f\|_{L_{p_2}(0,T;L_{p_1}(\Omega))}\cr 
&\quad+\|\omega(0)\|_{B_{p_1,p_2}^{2-2/p_2}(\Omega)}).\cr}
\leqno(2.6)
$$

\Proof 
To prove the lemma we use the idea of regularizer from \cite[Sect. 3] {17}, 
where all estimates are made in the H\"older spaces. Performing the estimates 
in the Sobolev spaces with the mixed norm (see \cite{18--21}) we prove the 
lemma.
\kwadrat

From \cite{22} we have

\proclaim Lemma 2.4. 
\item{(i)} Let $u\in W_{p,p_0}^{s,s/2}(\Omega^T)$, $s\in\R_+$, $s>2/p_0$, $p,p_0\in(1,\infty)$. 
Then $u(x,t_0)=u(x,t)|_{t=t_0}$ for $t_0\in[0,T]$ belongs to $B_{p,p_0}^{s-2/p_0}(\Omega)$ and
$$
\|u(\cdot,t_0)\|_{B_{p,p_0}^{s-2/p_0}(\Omega)}\le c\|u\|_{W_{p,p_0}^{s,s/2}(\Omega^T)},
$$
where constant $c$ does not depend on $u$.
\item{(ii)} For a given $\tilde u\in B_{p,p_0}^{s-2/p_0}(\Omega)$, $s\in\R_+$, $s>2/p_0$, $p,p_0\in(1,\infty)$, there exists a function $u\in W_{p,p_0}^{s,s/2}(\Omega^T)$ such that $u|_{t=t_0}=\tilde u$ for $t_0\in[0,T]$ and
$$
\|u\|_{W_{p,p_0}^{s,s/2}(\Omega^T)}\le c\|\tilde u\|_{B_{p,p_0}^{s-2/p_0}(\Omega)},
$$
where constant $c$ does not depend on $u$.

\section{3. Two-dimensional solutions}

First we have

\proclaim Lemma 3.1. 
Let $T>0$ be given. Assume that
$$
A_1^2\equiv{1\over\nu c_{s1}}\sup_{k\in\N_0}\intop_{kT}^{(k+1)T}
\|\bar f_s(t)\|_{L_2}^2dt<\infty,
\leqno\quad 1.
$$
$$
A_2^2\equiv{A_1^2\over1-e^{-\nu c_{s1}T}}+\|\bar v_s(0)\|_{L_2}^2<\infty,
\leqno\quad 2.
$$
where $c_{s1}$ is introduced in (2.3). Then
$$
\|\bar v_s(kT)\|_{L_2}^2\le A_2^2
\leqno(3.1)
$$
and
$$
\|\bar v_s(t)\|_{L_2}^2+\nu c_{s1}\intop_{kT}^t\|\bar v_s(t')\|_{H^1}^2dt'\le A_1^2+A_2^2\equiv A_3^2,
\leqno(3.2)
$$
where $t\in(kT,(k+1)T]$.
\goodbreak

\Proof 
Multiplying $(1.8)_1$ by $\bar v_s$, integrating over $\Omega$, using the 
periodic boundary conditions, the Poincar\'e inequality (2.3) and applying the 
Young inequality to the r.h.s. yield
$$
{d\over dt}\|\bar v_s\|_{L_2}^2+\nu c_{s1}\|\bar v_s\|_{H^1}^2\le{1\over\nu c_{s1}}
\|\bar f_s\|_{L_2}^2.
\leqno(3.3)
$$
Continuing, we obtain
$$
{d\over dt}(\|\bar v_s\|_{L_2}^2e^{\nu c_{s1}t})\le{1\over\nu c_{s1}}
\|\bar f_s\|_{L_2}^2e^{c_{s1}t}.
$$
Integrating with respect to time from $kT$ to $t\in(kT,(k+1)T]$ implies
$$
\|\bar v_s(t)\|_{L_2}^2\le{1\over\nu c_{s1}}\intop_{kT}^t\|\bar f_s(t')\|_{L_2}^2dt'+
e^{-\nu c_{s1}(t-kT)}\|\bar v_s(kT)\|_{L_2}^2.
$$
Setting $t=(k+1)T$ we get
$$
\|\bar v_s((k+1)T)\|_{L_2}^2\le{1\over\nu c_{s1}}\intop_{kT}^{(k+1)T}
\|\bar f_s(t)\|_{L_2}^2dt+e^{-\nu c_{s1}T}\|\bar v_s(kT)\|_{L_2}^2.
$$
By iteration we have
$$
\|\bar v_s(kT)\|_{L_2}^2\le{A_1^2\over1-e^{-\nu c_{s1}T}}+e^{-\nu c_{s1}kT}
\|v_s(0)\|_{L_2}^2\le A_2^2.
$$
Hence (3.1) is proved. Integrating (3.3) with respect to time from $t=kT$ to 
$t\in(kT,(k+1)T]$ and employing (3.1), we obtain (3.2). This concludes the 
proof.
\kwadrat

\noindent 
Next we obtain estimate for the second derivatives

\proclaim Lemma 3.2. 
Let assuptions of Lemma 3.1 hold. Let $\bar v_s(0)\in H^1(\Omega)$. Then
$$
\|\bar v_{sx}(kT)\|_{L_2}^2\le{c_{s1}A_1^2\over1-e^{-\nu c_{s1}T}}+
\|\bar v_{sx}(0)\|_{L_2}^2\equiv A_4^2
\leqno(3.4)
$$
and
$$
\|\bar v_{sx}(t)\|_{L_2}^2+\nu c_{s1}\intop_{kT}^t\|\bar v_s(t')\|_{H^2}^2dt'\le A_1^2+A_4^2\equiv A_5^2,
\leqno(3.5)
$$
where $t\in(kT,(k+1)T]$.
\goodbreak

\Proof 
Multiplying $(1.8)_1$ by $-\Delta\bar v_s$, integrating over $\Omega$ and using that $\bar v_s$ is divergence free yields
$$\eqal{
&-\intop_\Omega\bar v_{st}\cdot\Delta\bar v_sdx+\nu\intop_\Omega|\Delta\bar v_s|^2=\intop_\Omega v_s\cdot\nabla\bar v_s\cdot\Delta\bar v_sdx\cr 
&\quad-\intop_\Omega\bar f_s\cdot\Delta\bar v_sdx.\cr}
\leqno(3.6)
$$
Integrating by parts the first term on the l.h.s. equals
$$
{1\over2}{d\over dt}\intop_\Omega|\nabla\bar v_s|^2dx.
$$
To examine the first term on the r.h.s. of (3.6) we use the formula
$$
\Delta\bar v_s=\bigg(\matrix{&-(\rot\bar v_s)_{,x_2}\cr &(\rot\bar v_s)_{,x_1}\cr}\bigg)
$$
where $\rot\bar v_s=\bar v_{s2,x_1}-\bar v_{s1,x_2}$. Then
$$
\intop_\Omega v_s\cdot\nabla\bar v_s\cdot\Delta\bar v_sdx=\intop_\Omega(v_s\cdot
\nabla\bar v_{s2}\rot\bar v_{s,x_1}-v_s\cdot\nabla\bar v_{s1}\rot\bar v_{s,x_2})dx\equiv I.
$$
Performing integration by parts yields
$$\eqal{
I&=-\intop_\Omega(v_s\cdot\nabla\bar v_{s2,x_1}-v_s\cdot\nabla
\bar v_{s1,x_2})\rot\bar v_sdx\cr 
&\quad-\intop_\Omega(\bar v_{s,x_1}\cdot\nabla\bar v_{s2}-\bar v_{s,x_2}\cdot
\nabla\bar v_{s1})\rot\bar v_sdx\equiv I_1+I_2,\cr}
$$
where
$$
I_1=-\intop_\Omega v_s\cdot\nabla\rot\bar v_s\rot\bar v_sdx=0
$$
and
$$\eqal{
I_2&=-\intop_\Omega(\bar v_{s1,x_1}\bar v_{s2,x_1}+\bar v_{s2,x_1}
\bar v_{s2,x_2}-\bar v_{s1,x_2}\bar v_{s1,x_1}-\bar v_{s2,x_2}
\bar v_{s1,x_2})\rot\bar v_sdx\cr 
&=-\intop_\Omega\divv\bar v_s|\rot\bar v_s|^2dx=0.\cr}
$$
In view of the above considerations and the H\"older and Young inequalites applied to the last term on the r.h.s. of (3.6), we obtain from (3.6) the relation
$$
{d\over dt}\|\bar v_{sx}\|_{L_2}^2+\nu\|\Delta\bar v_s\|_{L_2}^2\le{1\over\nu}
\|\bar f_s\|_{L_2}^2.
\leqno(3.7)
$$
Applying the Poincar\'e inequality (see (2.3)) yields
$$
{d\over dt}\|\bar v_{sx}\|_{L_2}^2+\nu c_{s1}\|\bar v_{sx}\|_{L_2}^2\le
{1\over\nu}\|\bar f_s\|_{L_2}^2
$$
Hence
$$
{d\over dt}(\|\bar v_{sx}\|_{L_2}^2e^{\nu c_{s1}t})\le{1\over\nu}
\|\bar f_s\|_{L_2}^2e^{\nu c_{s1}t}
$$
Integrating with respect to time from $kT$ to $(k+1)T$ implies
$$
\|\bar v_{sx}((k+1)T)\|_{L_2}^2\le{1\over\nu}\intop_{kT}^{(k+1)T}
\|\bar f_s\|_{L_2}^2dt+e^{-\nu c_{s1}T}\|v_{sx}(kT)\|_{L_2}^2.
$$
Then iteration implies
$$
\|\bar v_{sx}(kT)\|_{L_2}^2\le{c_{s1}A_1^2\over1-e^{-\nu c_{s1}T}}+
e^{-\nu c_{s1}T}\|\bar v_{sx}(0)\|_{L_2}^2\le A_4^2
$$
Hence (3.4) is proved. Integrating (3.7) with respect to time from $kT$ to $t\in(kT,(k+1)T]$ and using Lemma 2.3 yields (3.5). This concludes the proof.
\kwadrat

\noindent 
To show stability of the two-dimensional solutions we need higher regularity of these solutions than the one proved in Lemma 3.2. We need such regularity that 
$v_s\in C(\R_+;W_\sigma^1(\Omega))$, $\sigma>3$. Moreover, we want to show that
$$
\|v_s(t)\|_{W_\sigma^1(\Omega)}\le c,
\leqno(3.8)
$$
where $c$ is a constant independent of time.

\noindent 
Finally, we do not want to apply the energy type method for higher derivatives 
(see \cite{14}), because it implied stronger restrictions on the external 
force. Instead, we are going to apply the increasing regularity technique. 
This is possible because in view of Lemma 3.2 the term 
$\bar v_s\cdot\nabla\bar v_s\in L_2(\R_+;L_\sigma(\Omega))$ with 
$\sigma\in(1,\infty)$. As it will be seen in Section 4 we will need to show 
only that $\bar v_s\in C(\R_+;W_\sigma^1(\Omega))$ with $\sigma>3$ (see (4.17)).

\proclaim Lemma 3.3. 
Assume that $v_s(0)\in B_{\sigma,2}^1(\Omega)$, 
$\bar f_s\in L_2(kT,(k+1)T;L_\sigma(\Omega))$, $k\in\N_0$, $\sigma>3$. Then 
$v_s\in C(\R_+;W_\sigma^1(\Omega))$, $\sigma>3$ and (3.8) holds.

\Proof 
Since $v_s\cdot\nabla\bar v_s\in L_2(kT,(k+1)T;L_\sigma(\Omega))$, \ 
$\sigma\in(1,\infty)$,\break $\bar f_s\in L_2(kT,(k+1)T;L_\sigma(\Omega))$, 
$\bar v_s(0)\in B_{\sigma,2}^1(\Omega)$ the theory from \cite{17--21},
and Lemma 3.2 imply the existence of solutions to (1.8) such that 
$v_s\in W_{\sigma,2}^{2,1}(\Omega\times\R_+)$ and the estimate holds
$$
\|\bar v_s\|_{W_{\sigma,2}^{2,1}(\Omega\times(0,T))}\le c(A_5^2+
\|\bar f_s\|_{L_2(0,T;L_\sigma(\Omega))}+\|v_s(0)\|_{B_{\sigma,2}^1(\Omega)}).
\leqno(3.9)
$$
However, we do not know how the constant $c$ depends on time. Therefore, we 
are not able to claim that (3.8) holds. Hence, we have to prove (3.8) step by 
step in time. Let us consider interval $(kT,(k+1)T)$. Let $\zeta=\zeta(t)$ be 
a smooth cut-off function such that $\zeta(t)=0$ for $t\in[kT,kT+\delta/2]$ 
and $\zeta(t)=1$ for $t\ge kT+\delta$. Introducing the new functions
$$
\tilde v_s=\bar v_s\zeta,\quad \tilde p_s=\bar p_s\zeta,\quad \tilde f_s=\bar f_s\zeta,\quad \dot\zeta=\zeta_{,t}
$$
we see that $(\tilde v_s,\tilde p_s)$ is a solution to the problem
$$\eqal{
&\tilde v_{st}-\nu\Delta\tilde v_s+\nabla\tilde p_s=\bar v_s\dot\zeta-v_s\cdot\nabla\tilde v_s+
\tilde f_s\quad &{\rm in}\ \ \Omega\times(kT,(k+1)T),\cr 
&\divv\tilde v_s=0\quad &{\rm in}\ \ \Omega\times(kT,(k+1)T),\cr 
&\tilde v_s|_{t=kT}=0.\cr}
\leqno(3.10)
$$
In view of Lemma 3.2 and \cite{17--21} we have the existence of solutions to 
(3.10) such that
$$
\tilde v_s\in W_{\sigma,2}^{2,1}(kT+\delta,(k+1)T;\Omega),\quad 
\nabla\tilde p_s\in L_{\sigma,2}(kT+\delta,(k+1)T;\Omega)
$$
and the estimate holds
$$
\|\tilde v_s\|_{W_{\sigma,2}^{2,1}(kT+\delta,(k+1)T;\Omega)}\le c
\bigg({1\over\delta}A_5+A_5^2+\|\tilde f_s\|_{L_2(kT+\delta/2,(k+1)T;\Omega)}),
\leqno(3.11)
$$
where $c$ might depend on $T$ but it does not depend on $k$. Hence by 
imbedding for $\sigma>3$ estimate (3.11) implies (3.8).

To get (3.8) we need only estimate for the interval $(kT,kT+\delta)$, 
$k\in\N$, because for $k=0$ we have (3.9). From (3.11) for $k$ replaced by 
$k-1$ we obtain the estimate for
$$
\|\tilde v_s\|_{W_{\sigma,2}^{2,1}((k-1)T+\delta,kT;\Omega)}\le c\bigg({1\over\delta}A_5+A_5^2+
\|\tilde f_s\|_{L_2((k-1)^+{\delta\over2},kT;\Omega)}\bigg),
$$
so by the trace thoerem (see Lemma 2.3) we derive
$$
\|\tilde v_s(kT)\|_{B_{\sigma,2}^1(\Omega)}\le c
\|\bar v_s\|_{W_{\sigma,2}^{2,1}((k-1)T+\delta,kT,\Omega)}.
\leqno(3.12)
$$
Hence, repeating the considerations leading to (3.9) for time interval\\ 
$(kT,kT+\delta)$ we obtain that $\tilde v_s\in W_{\sigma,2}^{2,1}(\Omega\times
(kT,kT+\delta))$ and the estimate
$$\eqal{
&\|\tilde v_s\|_{W_{\sigma,2}^{2,1}(\Omega\times(kT,kT+\delta))}\le c(A_5^2+
\|\bar f_s\|_{L_{\sigma,2}(\Omega\times(kT,kT+\delta))}\cr 
&\quad+\|\bar v_s(kT)\|_{B_{\sigma,2}^1(\Omega)}).\cr}
\leqno(3.13)
$$
Hence (3.8) holds for all $t\in\R_+$ and Lemma 3.3 is proved.
\kwadrat

\section{4. Stability}

In this Section we examine problem (1.4). First we derive a global estimate 
for $L_2$ norm of $u$. We show how appears restriction from Assumption 2 of 
Lemma 4.1 (it is much more restrictive in \cite{14}). Fortunately, we do not 
need Lemma 4.1 to prove stability. Hence we have

\proclaim Lemma 4.1. 
Let the assumptions of Lemmas 3.1, 3.2 hold. Let
\item{1.} $B_1^2=\sup_{k\in\N_0}\intop_{kT}^{(k+1)T}
\big({\nu c_1\over2c_3}\big|\intop_0^t
\diagintop_\Omega g(t')dxdt'+\diagintop_\Omega u(0)dx\big|^2\\
\hbox to .6cm{}+{2c_3\over\nu c_1}
\|\bar g(t)\|_{L_{6/5}}^2\big)dt$, where $c_1$ follows from Poincar\'e 
inequality (2.4) and $c_3$ from imbedding (4.3).
\item{2.} ${-\nu c_1\over2}T+{4c_3\over\nu c_1}A_3^2\le0$.
\vskip6pt
\item{3.} $B_2^2=\exp\big({4c_3\over\nu c_1}A_5^2\big)B_1^2$.\\
Then
$$\eqal{
&\|\bar u(kT)\|_{L_2}^2\le{B_2^2\over1-\exp(-\nu c_1T/2)}+
\|\bar u(0)\|_{L_2}^2=B_3^2,\cr 
&\|\bar u(t)\|_{L_2}^2\le B_2^2+B_3^2\equiv B_4^2,\cr}
\leqno(4.1)
$$
for $t\in[kT,(k+1)T]$ and any $k\in\N_0$.

\Proof 
Multiplying $(1.9)_1$ by $\bar u$ and integrating over $\Omega$ gives
$$\eqal{
&{1\over2}{d\over dt}\|\bar u\|_{L_2}^2+\nu c_1\|\bar u_x\|_{L_2}^2\le
\bigg|\intop_\Omega u\cdot\nabla\bar v_s\cdot\bar udx\bigg|+
\bigg|\intop_\Omega\bar g\cdot\bar udx\bigg|\cr 
&\le\bigg|\intop_\Omega\bar u\cdot\nabla\bar v_s\cdot\bar udx\bigg|+
\bigg|\diagintop_\Omega udx\cdot\intop_\Omega\nabla\bar v_s\cdot\bar udx\bigg|
+\bigg|\intop_\Omega\bar g\cdot\bar udx\bigg|.\cr}
\leqno(4.2)
$$
Employing the estimates
$$\eqal{
&\bigg|\intop_\Omega\bar u\cdot\nabla\bar v_s\cdot\bar udx\bigg|\le{\varepsilon_1\over2}
\|\bar u\|_{L_6}^2+{1\over2\varepsilon_1}\|\bar v_{sx}\|_{L_3}^2\|\bar u\|_{L_2}^2,\cr 
&\bigg|\diagintop udx\cdot\intop_\Omega\nabla\bar v_s\cdot\bar udx\bigg|={1\over2\varepsilon_2}
\|\nabla\bar v_s\|_{L_2}^2\|\bar u\|_{L_2}^2+{\varepsilon_2\over2}
\bigg|\diagintop udx\bigg|^2,\cr 
&\bigg|\intop_\omega\bar g\cdot\bar udx\bigg|\le{\varepsilon_3\over2}\|\bar u\|_{L_6}^2+
{1\over2\varepsilon_3}\|\bar g\|_{L_{6/5}}^2,\cr}
$$
(2.4), the imbedding
$$
\|\bar u\|_{L_6}^2\le c_3\|\bar u\|_{H^1}^2
\leqno(4.3)
$$
and that ${\varepsilon_i\over2}$ $c_3\le{c_1\over4}$, $i=1,3$, we obtain from (4.2) the inequality
$$\eqal{
&{1\over2}{d\over dt}\|\bar u\|_{L_2}^2+{\nu c_1\over2}\|\bar u\|_{H^1}^2\le{c_3\over\nu c_1}
(\|\bar v_{sx}\|_{L_3}^2+\|\bar v_{sx}\|_{L_2}^2)\|\bar u\|_{L_2}^2\cr 
&\quad+{\nu c_1\over4c_3}\bigg|\diagintop udx\bigg|^2+{c_3\over\nu c_1}
\|\bar g\|_{L_{6/5}}^2,\cr}
\leqno(4.4)
$$
where $\varepsilon_2={\nu c_1\over2c_3}$ is set. Employing (2.2) in (4.4) yields
$$\eqal{
&{d\over dt}\|\bar u\|_{L_2}^2+\nu c_1\|\bar u\|_{H^1}^2\le{2c_3\over\nu c_1}
(\|\bar v_{sx}\|_{L_3}^2+\|\bar v_{sx}\|_{L_2}^2)\|\bar u\|_{L_2}^2\cr 
&\quad+{\nu c_1\over2c_3}\bigg|\intop_0^t\diagintop_\Omega g(t')dxdt'+
\diagintop_\Omega u(0)dx\bigg|^2+{2c_3\over\nu c_1}\|\bar g\|_{L_{6/5}}^2.\cr}
\leqno(4.5)
$$
Considering (4.5) for $t\in(kT,(k+1)T)$ we have
$$\eqal{
&{d\over dt}\bigg[\|\bar u(t)\|_{L_2}^2\exp\bigg(\nu c_1t-{2c_3\over\nu c_1}\intop_{kT}^t
(\|\bar v_{sx}(t')\|_{L_3}^2+\|\bar v_{sx}(t')\|_{L_2}^2)dt'\bigg)\bigg]\cr 
&\le\bigg({\nu c_1\over2c_3}\bigg|\intop_0^t\diagintop_\Omega g(t')dxdt'+
\diagintop_\Omega u(0)dx\bigg|^2+{2c_3\over\nu c_1}\|\bar g(t)\|_{L_{6/5}}^2\bigg)\cdot\cr 
&\quad\cdot\exp\bigg(\nu c_1t-{2c_3\over\nu c_1}\intop_{kT}^t
(\|\bar v_{sx}(t')\|_{L_3}^2+\|\bar v_{sx}(t')\|_{L_2}^2)dt'\bigg).\cr}
\leqno(4.6)
$$
Integrating (4.6) with respect to time from $t=kT$ to $t\in(kT,(k+1)T]$ implies
$$\eqal{
&\|\bar u(t)\|_{L_2}^2\le\exp\bigg[{2c_3\over\nu c_1}\intop_{kT}^t
(\|v_{sx}(t')\|_{L_3}^2+\|v_{sx}(t')\|_{L_2}^2)dt'\bigg]\cdot\cr 
&\quad\cdot\intop_{kT}^t\bigg({\nu c_1\over2c_3}\bigg|\intop_0^{t'}\diagintop_\Omega g(t'')dxdt''+\diagintop_\Omega u(0)dx\bigg|^2+{2c_3\over\nu c_1}
\|\bar g(t')\|_{L_{6/5}}^2\bigg)dt'\cr 
&\quad+\|\bar u(kT)\|_{L_2}^2\exp\bigg[-\nu c_1(t-kT)+{2c_3\over\nu c_1}\intop_{kT}^t
(\|\bar v_{sx}(t')\|_{L_3}^2\cr 
&\quad+\|\bar v_{sx}(t')\|_{L_2}^2)dt'\bigg].\cr}
\leqno(4.7)
$$
Setting $t=(k+1)T$ and using (3.5), inequality (4.7) yields
$$\eqal{
&\|\bar u((k+1)T)\|_{L_2}^2\cr 
&\le\exp\bigg({4c_3\over\nu c_1}A_3^2\bigg)\intop_{kT}^{(k+1)T}
\bigg[{\nu c_1\over2c_3}\bigg|\intop_0^t\diagintop_\Omega g(t')dxdt'+
\diagintop_\Omega u(0)dx\bigg|^2\cr 
&\quad+{2c_3\over c_1}\|\bar g(t)\|_{L_{6/5}}^2\bigg]dt+\|\bar u(kT)\|_{L_2}^2\exp 
\bigg(-\nu c_1T+{4c_3\over\nu c_1}A_3^2\bigg).\cr}
\leqno(4.8)
$$
In view of assumptions 1--3 of the lemma we have
$$
\|\bar u((k+1)T)\|_{L_2}^2\le B_2^2+\exp\bigg({-\nu c_1\over2}T\bigg)\|\bar u(kT)\|_{L_2}^2.
\leqno(4.9)
$$
Iteration implies
$$
\|\bar u(kT)\|_{L_2}^2\le{B_2^2\over1-\exp(-\nu c_1T/2)}+\exp\bigg({-\nu c_1\over2}kT\bigg)
\|\bar u(0)\|_{L_2}^2.
\leqno(4.10)
$$
Hence $(4.1)_1$ is proved. Employing assumptions of the lemma and $(4.1)_1$ 
in (4.7) gives $(4.1)_2$. This concludes the proof.
\kwadrat

\Remark{4.2.} 
Assumption 2 of Lemma 4.1 has the explicit form
$$\eqal{
&{2-\exp(-\nu c_{s1}T)\over c_{s1}\nu(1-\exp(-\nu c_{s1}T))}\sup_{k\in\N_0}\intop_{kT}^{(k+1)T}
\|\bar f_s(t)\|_{L_2}^2dt\cr 
&\quad+\|\bar v_{sx}(0)\|_{L_2}^2\le{\nu^2c_1^2\over8c_3}T.\cr}
\leqno(4.11)
$$
Assuming that $\|\bar v_{sx}(0)\|_{L_2}$ is given we see that (4.11)holds for\\ 
$T>{8c_3\over\nu^2c_1^2}\|\bar v_{sx}(0)\|_{L_2}^2$.
For such large $T$ we have a strong restriction on 
$\sup_{k\in\N_0}\intop_{kT}^{(k+1)T}\|\bar f_s(t)\|_{L_2}^2dt$. Physically, it means that the energy introduced to the considered region must not to be too large comparing with the dissipation.

Finally, we show that 3d solutions to (1.1) remain close to 2d solutions to 
(1.2) for all time if their initial data and the external forces are 
sufficiently close. In this proof we omit the heavy restriction (4.11).

\proclaim Lemma 4.3. 
Let $\bar v_s\in C(\R_+;W_3^1)$, $\bar g\in C(\R_+;L_2)$, 
$\bar u(0)\in H^1$. Let $\gamma\in(0,\gamma_*]$, where 
$\nu c_4-{c_5\over\nu^3}\gamma_*^2\ge{c_*\over2}$, $c_*<\nu c_4$ and $c_4$, 
$c_5$ are introduced in (4.16). Assume that
$$\eqal{
&\|\bar u(0)\|_{H^1}^2\le\gamma\cr 
&G^2(t)={c_5\over\nu}\bigg[\|\bar v_{sx}\|_{L_3}^2
\bigg|\intop_0^t\diagintop_\Omega g(x,t')dxdt'+\diagintop_\Omega u(0)dx\bigg|^2\cr 
&\quad+\|\bar g\|_{L_2}^2\bigg]\le c_*{\gamma\over4}.\cr}
\leqno(4.12)
$$
Let $T>0$ be given and $k\in\N_0$. assume that
$$\eqal{
&{c_5\over\nu}\intop_{kT}^{(k+1)T}\|\bar v_{sx}\|_{L_3}^2dt\le{c_*\over4}T,\qquad \intop_{kT}^{(k+1)T}G^2(t)dt\le\alpha\gamma,\cr
&\alpha\exp\bigg({c_*\over4}T\bigg)+\exp\bigg(-{c_*\over4}T\bigg)\le 1.\cr}
$$
Then
$$
\|\bar u(t)\|_{H^1}^2\le\gamma\quad {\sl for}\ \ t\in\R_+.
\leqno(4.13)
$$

\Proof 
Differentiating $(1.9)_1$ with respect to $x$, multiplying the result by $\bar u_x$, integrating over $\Omega$ and employing the periodic boundary conditions yield
$$\eqal{
&{1\over2}{d\over dt}\|\bar u_x\|_{L_2}^2+\nu\|\bar u_{xx}\|_{L_2}^2\le\|\bar u_x\|_{L_3}^3+
\bigg|\intop_\Omega\bar v_{sx}\cdot\nabla u\cdot\bar u_xdx\bigg|\cr 
&\quad+2\bigg|\intop_\Omega\bar u_x\cdot\nabla\bar v_s\cdot\bar u_xdx\bigg|+
\bigg|\intop_\Omega u\cdot\nabla\bar v_s\cdot\bar u_{xx}dx\bigg|+
\bigg|\intop_\Omega\bar g\cdot\bar u_{xx}dx\bigg|.\cr}
\leqno(4.14)
$$
Adding (4.2) and (4.14), applying the H\"older, the Young and the Poincar\'e 
inequalities, we derive
$$\eqal{
&{d\over dt}\|\bar u\|_{H^1}^2+\nu c\|\bar u\|_{H^2}^2\le c(\|\bar u_x\|_{L_3}^3+
{1\over\nu}\|\bar v_{sx}\|_{L_3}^2\|\bar u_x\|_{L_2}^2\cr 
&\quad+{1\over\nu}\|u\|_{L_6}^2\|\bar v_{sx}\|_{L_3}^2+{1\over\nu}
\|\bar g\|_{L_2}^2).\cr}
$$
Using $\|u\|_{L_6}^2\le c\big(\|\bar u\|_{L_6}^2+
\big|\diagintop_\Omega udx\big|^2\big)$ and $\|\bar u\|_{L_6}\le 
c\|\bar u\|_{H^1}\le c\|\bar u_x\|_{L_2}$, which holds in view of the 
Poincar\'e inequality, we get
$$\eqal{
&{d\over dt}\|\bar u\|_{H^1}^2+\nu c\|\bar u\|_{H^2}^2\cr 
&\le c\bigg[\|\bar u_x\|_{L_3}^3+{1\over\nu}
\|\bar v_{sx}\|_{L_3}^2\bigg(\|\bar u_x\|_{L_2}^2+\bigg|\diagintop_\Omega udx\bigg|^2\bigg)
+{1\over\nu}\|\bar g\|_{L_2}^2\bigg].\cr}
\leqno(4.15)
$$
In view of (2.2) and the interpolation inequality (see \cite[Ch. 3, Sect. 15] 
{16})
$$
\|\bar u_x\|_{L_3}\le c\|\bar u_{xx}\|_{L_2}^{1/2}\|\bar u_x\|_{L_2}^{1/2}
$$
(which holds without the lower order term because 
$\intop_\Omega\bar u_xdx=0$), we obtain from (4.15) the inequality
$$\eqal{
&{d\over dt}\|\bar u\|_{H^1}^2+\nu c_4\|\bar u\|_{H^2}^2\le{c_5\over\nu^3}
\|\bar u_x\|_{L_2}^6+{c_5\over\nu}\|\bar v_{sx}\|_{L_3}^2\|\bar u_x\|_{L_2}^2\cr 
&\quad+{c_5\over\nu}\|\bar v_{sx}\|_{L_3}^2\bigg|\intop_0^t\diagintop_\Omega 
g(x,t')dxdt'+\diagintop_\Omega u(0)dx\bigg|^2+
{c_5\over\nu}\|\bar g\|_{L_2}^2.\cr}
\leqno(4.16)
$$
To prove the lemma we need to know that the r.h.s. of (4.16) is bounded. 
We consider inequality (4.16) in the time interval $(kT,(k+1)T)$, $k\in\N_0$. 
Assume that we have proved that $u(kT)\in H^1(\Omega)$ and 
$\|u(kT)\|_{H^1}^2\le\gamma$, where $\gamma$ is sufficiently small. Using that 
$g\in L_2(\Omega\times(kT,(k+1)T))$ is sufficiently small we have existence of 
solutions to problem (1.9) in $W_2^{2,1}(\Omega\times(kT,(k+1)T))$ because the 
other terms on the r.h.s. of (1.9) also belong to 
$L_2(\Omega\times(kT,(k+1)T))$ in view of imbeddings and assumption that 
$v_s\in W_2^{2,1}(\Omega\times(kT,(k+1)T))$. The last assertion holds in view 
of the assumptions of Lemma 3.2 and the restriction that $v_s$ is 
a two-dimensional solution to the Navier-Stokes equations. However, to have 
the r.h.s. of (4.16) bounded we need that 
$v_s\in L_\infty(kT,(k+1)T;W_{3^+}^1(\Omega))$, where $3^+>3$ but close to 3. 
This follows from Lemma 3.3, where it is proved that 
$v_s\in W_{\sigma,2}^{2,1}(\Omega\times(kT,(k+1)T))$ for any $\sigma$ if data 
are sufficiently smooth.

\noindent 
In view of the above remarks we can introduce the quantities
$$\eqal{
G^2(t)&={c_5\over\nu}\bigg(\|\bar v_{sx}\|_{L_3}^2
\bigg|\intop_0^t\diagintop_\Omega g(x,t')dxdt'+\diagintop_\Omega u(0)dx\bigg|^2\cr 
&\quad+{1\over\ }\|\bar g\|_{L_2}^2\bigg),\quad 
A^2(t)={c_5\over\nu}\|\bar v_{sx}\|_{L_3}^2,\cr 
&X(t)=\|\bar u(t)\|_{H^1},\quad Y(t)=\|\bar u(t)\|_{H^2}.\cr}
\leqno(4.17)
$$
Then (4.17) takes the form
$$
{d\over dt}X^2+\nu c_4Y^2\le {c_5\over\nu^3}X^4X^2+A^2X^2+G^2.
$$
Since $X\le Y$ we have
$$
{d\over dt}X^2\le-X^2\bigg(\nu c_4-{c_5\over\nu^3} X^4\bigg)+A^2X^2+G^2
\leqno(4.18)
$$
Let $\gamma\in(0,\gamma_*]$, where $\gamma_*$ is so small that
$$
\nu c_4-{c_5\over\nu^3}\gamma_*^2\ge c_*/2,\quad c_*<\nu c_4.
\leqno(4.19)
$$
Since the coefficients of equation (4.18) depend on the two-dimensional 
solution determined step by step in time we consider (4.18) in the interval 
$[kT,(k+1)T]$, $k\in\N_0$, with the assumptions
$$
X^2(kT)\le\gamma,\quad G^2(t)\le c_*\gamma/4\quad {\rm for\ all}\ \ 
t\in[kT,(k+1)T].
$$
Let us introduce the quantity
$$
Z^2(t)=\exp\bigg(-\intop_{kT}^tA^2(t')dt'\bigg)X^2(t),\quad t\in[kT,(k+1)T].
$$
Then (4.18) takes the form
$$
{d\over dt}Z^2\le-\bigg(\nu c_4-{c_5\over\nu^3}X^4\bigg)Z^2+\bar G^2,
\leqno(4.20)
$$
where $\bar G^2=G^2\exp\big(-\intop_{kT}^tA^2(t')dt'\big)$.

\noindent
Suppose that
$$\eqal{
t_*&=\inf\{t\in(kT,(k+1)T]:\>X^2(t)>\gamma\}\cr 
&=\inf\bigg\{t\in(kT,(k+1)T]:\>Z^2(t)
>\gamma\exp\bigg(-\intop_{kT}^tA^2(t')dt'\bigg)\bigg\}>kT.\cr}
$$
By (4.19) for $t\in(0,t_*]$ inequality (4.20) takes the form
$$
{d\over dt}Z^2\le-{c_*\over2}Z^2+\bar G^2(t).
\leqno(4.21)
$$
Clearly, we have
$$\eqal{
&Z^2(t_*)=\gamma\exp\bigg(-\intop_{kT}^{t_*}A^2(t')dt'\bigg)\quad 
&{\rm and}\quad\cr 
&Z^2(t)>\gamma\exp\bigg(-\intop_{kT}^{t_*}A^2(t')dt'\bigg)\quad \cr
&{\rm for}\quad t>t_*.\cr}
\leqno(4.22)
$$
Then (4.21) yields
$$
{d\over dt}Z^2|_{t=t_*}\le c_*\bigg(-{\gamma\over2}+{\gamma\over4}\bigg)
\exp\bigg(-\intop_{kT}^{t_*}A^2(t')dt'\bigg)<0
$$
contradicting with (4.22). Therefore
$$
Z^2(t)<\gamma\exp\bigg(-\intop_{kT}^{t_*}A^2(t')dt'\bigg)\quad {\rm for}\quad 
t>t_*.
\leqno(4.23)
$$
Then definition of $Z^2(t)$ implies
$$
X^2(t)\le\gamma\exp\bigg(\intop_{t_*}^tA^2(t')dt'\bigg)\quad {\rm for}\quad t>t_*.
$$
For sufficiently small $\gamma$ inequality (4.18) takes the form
$$
{d\over dt}X^2+{c_*\over2}X^2\le A^2X^2+G^2.
\leqno(4.24)
$$
Integrating (4.24) with respect to time from $t=kT$ to $t=(k+1)T$ gives
$$\eqal{
X^2((k+1)T)&\le\exp\bigg(\intop_{kT}^{(k+1)T}A^2(f)dt\bigg)\intop_{kT}^{(k+1)T}
G^2(t)dt\cr
&\quad+\exp\bigg(-{c_*\over2}T+\intop_{kT}^{(k+1)T}A^2(t)dt\bigg)X^2(kT).\cr}
\leqno(4.25)
$$
In view of the assumptions
$$
{c_*\over4}T\ge\intop_{kT}^{(k+1)T}A^2(t)dt,\quad 
\intop_{kT}^{(k+1)T}G^2(t)dt\le\alpha\gamma,
\leqno(4.26)
$$
where $\alpha$ is so small and $T$ so targe that
$$
\alpha\exp\bigg(\intop_{kT}^{(k+1)T}A^2(t)dt\bigg)+\exp\bigg(-{c_*\over4}T\bigg) \le 1,
\leqno(4.27)
$$
we have that $X^2((k+1)T)<\gamma$. Then by the induction we prove the lemma.\hfill
{\kwadrat}

\subsection{Acknowledgements}

The author thanks to professor Y.Shibata for very important comments concerning the proof of Lemma 4.3.
\vskip6pt
The research leading to these results has received funding from the People Programme (Marie Curie Actions) of the European Union's Seventh Framework Programme FP7/2007-2013/ under REA grant agreement n$^\circ$~319012 and from the Funds for International Co-operation under Polish Ministry of Science and Higher Education grant agreement n$^\circ$~2853/7.PR/2013/2.

\section{References}

\item{1.} Beir\~ao da Veiga H. and Secchi P.: $Lp$-stability for the strong 
solutions of the Navier-Stokes equations in the whole space, Arch. Ration. 
Mech. Anal. 98 (1987), 65--69.

\item{2.} Ponce G., Racke R., Sideris T. C. and Titi E. S.: Global stability 
of large solutions to the 3d Navier-Stokes equations, Comm. Math. Phys. 159 
(1994), 329--341.

\item{3.} Mucha P. B.: Stability of 2d incompressible flows in $R^3$, J. 
Diff. Eqs. 245 (2008), 2355--2367.

\item{4.} Gallagher I.: The tridimensional Navier-Stokes equations with 
almost bidimensional data: stability, uniqueness and life span, Internat. Mat. 
Res. Notices 18 (1997), 919--935.

\item{5.} Iftimie D.: The 3d Navier-Stokes equations seen as a perturbation 
of the 2d Navier-Stokes equations, Bull. Soc. Math. France 127 (1999), 
473--517.

\item{6.} Mucha P. B.: Stability of nontrivial solutions of the Navier-Stokes 
system on the three-dimensional torus, J. Differential Equations 172 (2001), 
359--375.

\item{7.} Mucha P. B.: Stability of constant solutions to the Navier-Stokes 
system in $\R^3$, Appl. Math. 28 (2001), 301--310.

\item{8.} Auscher P., Dubois S. and Tchamitchian P.: On the stability of 
global solutions to Navier-Stokes equations in the space, Journal de 
Math\'ematiques Pures et Appliques 83 (2004), 673--697.

\item{9.} Gallagher I., Iftimie D. and Planchon F.: Asymptotics and stability 
for global solutions to the Navier-Stokes equations, Annales de I'Institut 
Fourier 53 (2003), 1387--1424.

\item{10.} Zhou Y.: Asymptotic stability for the Navier-Stokes equations 
in the marginal class, Proc Roy. Soc. Edinburgh 136 (2006), 1099--1109.

\item{11.} Karch G. and Pilarczyk D.: Asymptotic stability of Landau 
solutions to Navier-Stokes system, Arch. Rational Mech. Anal. 202 (2011), 
115--131.

\item{12.} Chemin J. I, and Gallagher I.: Wellposedness and stability results 
for the Navier-Stokes equations in $R^3$, Ann. I. H. Poincar\'e, Analyse Non 
Lin\'eaire, 26 (2009), 599--624.

\item{13.} Bardos C., Lopes Filho M. C., Niu D., Nussenzveig Lopes H. J. and 
Titi E. S.: Stability of two-dimensional viscous incompressible flows under 
three-dimesional perturbations and inviscid symmetry breaking, SIAM J. Math. 
Anal. 45 (2013), 1871--1885.

\item{14.} Zadrzy\'nska, E.; Zaj\c aczkowski, W. M.: Stability of 
two-dimensional Navier-Stokes motions in the periodic case, J. Math. Anal. 
Appl., 423 (2015), 956--974.

\item{15.} Adams, R. A.; Fournier, J. J. F.: Sobolev spaces, Second Edition, 
Academic Press 2008.

\item{16.} Besov, O. V.; Il'in, V. P.; Nikol'skii, S. M.: Integral 
representations of functions and imbedding theorems, Nauka, Moscow 1975 
(in Russian).

\item{17.} Solonnikov, V. A.: On the solvability of generalized Stokes equations in the spaces of periodic functions, Ann. Univ. Ferrara-Sez. VII-Sc. Mat. XLVI (2000), 219--249.

\item{18.} Krylov, N. V.: The heat equation in $L_q(0,T;L_p)$-spaces with weights, SIAM J. Math. Anal. 32 (2001), 1117--1141.

\item{19.} Krylov, N. V.: The Calderon-Zygmund theorem and its application to parabolic equations, Algebra i Analiz 13 (2001), 1--25 (in Russian).

\item{20.} Solonnikov, V. A.: Estimates of solutions of the Stokes equations in Sobolev spaces with a mixed norm, Zap. Nauchn. Sem. POMI 288 (2002), 204--231.

\item{21.} Solonnikov, V. A.: On the estimates of solutions of nonstationary Stokes problem in anisotropic S.L. Sobolev spaces and on the estimate of resolvent of the Stokes problem, Usp. Mat. Nauk. 58 (2) (2003) (350), 123--156.

\item{22.} Bugrov, Y. S.: Function spaces with mixed norm, Math. USSR-Izv. 5 (1971), 1145--1167 (in Russian).

\bye